\documentclass[12pt,reqno]{amsart} 
\usepackage[T1]{fontenc}
\usepackage[utf8]{inputenc}
\usepackage{lmodern,microtype}
\usepackage{amsmath,amssymb,amsthm,mathtools}
\usepackage{bbm}
\usepackage[margin=1.25in]{geometry}
\usepackage{graphicx,booktabs}
\usepackage[dvipsnames]{xcolor}
\usepackage{tikz}

\usepackage[
  backend=biber,
  style=alphabetic,
  giveninits=true,
  maxnames=10,
  maxalphanames=10,
  doi=false,url=false,eprint=false,isbn=false
]{biblatex}
\DeclareSourcemap{
  \maps[datatype=bibtex]{
    \map{
      \step[fieldsource=note, match=\regexp{\\url}, final]
      \step[fieldsource=note, fieldtarget=addendum]
    }
    \map{
      \step[fieldsource=howpublished, match=\regexp{\\arXiv}, final]
      \step[fieldsource=howpublished, fieldtarget=addendum]
    }
  }
}
\renewbibmacro*{addendum+pubstate}{%
  \setunit{\addperiod\space}%
  \printfield{addendum}%
  \newunit\newblock
  \printfield{pubstate}}

\DeclareNameAlias{default}{given-family}
\DeclareNameAlias{sortname}{given-family}
\DeclareFieldFormat[article,inbook,incollection,inproceedings,misc,unpublished]
  {title}{\mkbibemph{#1}}
\DeclareFieldFormat{journaltitle}{#1\isdot}
\DeclareFieldFormat{booktitle}{#1\isdot}
\DeclareFieldFormat{pages}{#1}
\DeclareFieldFormat[article]{volume}{\mkbibbold{#1}}
\DeclareFieldFormat[article]{number}{\bibstring{number}\addnbspace#1}

\renewbibmacro*{in:}{}
\renewbibmacro*{journal+issuetitle}{%
  \usebibmacro{journal}%
  \setunit*{\addspace}%
  \printfield{series}%
  \setunit*{\addspace}%
  \printfield{volume}%
  \setunit{\addspace}%
  \usebibmacro{issue+date}%
  \setunit{\addcomma\space}%
  \printfield{number}%
  \setunit{\addcomma\space}%
  \printfield{eid}%
  \newunit}
  
\usepackage[colorlinks=true,linkcolor=Red,citecolor=ForestGreen,urlcolor=ForestGreen]{hyperref}
\def\arXiv#1{\href{http://arxiv.org/abs/#1}{\tt arXiv:#1}}

\hypersetup{pdftitle={Asymptotic isoperimetry for parallelohedra}}

\newtheorem{theorem}{Theorem}
\newtheorem{proposition}[theorem]{Proposition}
\newtheorem{lemma}[theorem]{Lemma}
\newtheorem{corollary}[theorem]{Corollary}
\newtheorem{question}[theorem]{Question}

\theoremstyle{definition}

\newcommand{\R}{\mathbb R}
\newcommand{\E}{\mathbb E}
\newcommand{\Z}{\mathbb Z}
\newcommand{\Prb}{\mathbb P}
\newcommand{\PP}{\mathrm P}
\newcommand{\AV}{\mathrm A}
\newcommand{\VV}{\mathrm V}
\newcommand{\J}{\mathcal J}
\newcommand{\ind}{\mathbbm{1}}
\DeclareMathOperator{\iq}{iq}
\DeclareMathOperator{\area}{area}
\DeclareMathOperator{\vol}{vol}
\DeclareMathOperator{\Vor}{Vor}
\DeclareMathOperator{\Var}{Var}

\DeclareMathOperator{\argmax}{argmax}
\allowdisplaybreaks[2]
\title{Asymptotic isoperimetry for parallelohedra}

\author{Matteo Novaga}
\address{Department of Mathematics, University of Pisa, Pisa, Italy}
\email{matteo.novaga@unipi.it}

\author{Lark Song}
\address{Department of Mathematics, University of Wisconsin, Madison, WI, USA}
\email{lark.song@wisc.edu}

\begin{document}
\begin{abstract}
We prove that the smallest possible surface area of a convex lattice tile of $\R^n$ approaches $\sqrt{2}$ times that of a ball of equal volume as $n\to\infty$.
For the lower bound, we combine Gaussian smoothing with the Milman-Neeman multi-bubble theorem. For the upper bound, we use Siegel mean value theorem to average over lattice Voronoi cells.
Our method also yields convex $\Z^n$-tiles whose orthogonal projection onto any hyperplane has volume $O(\sqrt{\log n})$, answering a question of Naor and Regev.
\end{abstract}

\maketitle

\thispagestyle{empty}

\section{Introduction}\label{sec:introduction}

The isoperimetric problem asks what shape minimizes surface area at fixed volume.
The ball gives this minimum but cannot tile space.
We consider convex bodies that tile space by translations.
Minkowski proved that such a convex body must be a centrally symmetric polytope with centrally symmetric facets \cite{Minkowski1897}. Minkowski and Venkov proved that every convex translational tile admits a lattice tiling \cite{Minkowski1897,Venkov1954}, as also proved by McMullen \cite{McMullen1980,McMullen1981}. A convex translational tile is called a \emph{parallelohedron}, some literature uses the term \emph{parallelotope}.

\begin{figure}[!ht]
\centering
\begin{tikzpicture}[x=1cm,y=1cm,line join=round,line cap=round,
                    shape drawing/.style={scale=1.25},
                    every node/.style={inner sep=0pt},font=\small]
  \draw[shape drawing,line width=.6pt] (0:.72)
    \foreach \angle in {60,120,180,240,300} { -- (\angle:.72) } -- cycle;

  \begin{scope}[shift={(3.70,0)},shape drawing]
    \begin{scope}[scale=.75,x={(.380cm,-.178cm)},y={(-.325cm,-.207cm)},
                  z={(0cm,.419cm)},line width=.55pt]
      \draw (2,1,0) -- (2,0,1) -- (2,-1,0) -- (2,0,-1) -- cycle;
      \draw (1,2,0) -- (0,2,1) -- (-1,2,0) -- (0,2,-1) -- cycle;
      \draw (1,0,2) -- (0,1,2) -- (-1,0,2) -- (0,-1,2) -- cycle;
      \foreach \sx/\sy/\sz in {1/1/1,-1/1/1,1/-1/1,1/1/-1}{
        \draw ({2*\sx},\sy,0) -- (\sx,{2*\sy},0) -- (0,{2*\sy},\sz)
          -- (0,\sy,{2*\sz}) -- (\sx,0,{2*\sz}) -- ({2*\sx},0,\sz) -- cycle;
      }
    \end{scope}
  \end{scope}

  \begin{scope}[shift={(7.40,0)},shape drawing,
                rear wire/.style={black,line width=.22pt},
                wire/.style={black,line width=.32pt},
                near wire/.style={black,line width=.40pt,
                                  preaction={draw=white,line width=.85pt}}]
    \foreach \vertex/\u/\v in {
      outerxp/-0.906/0.599, outerxm/0.904/-0.332,
      outeryp/0.784/0.699, outerym/-0.489/-0.362,
      outerzp/-0.063/-0.860, outerzm/0.166/0.861,
      innerxp/-0.118/0.194, innerxm/0.242/0.009,
      inneryp/0.194/0.207, innerym/-0.057/-0.002,
      innerzp/0.042/-0.074, innerzm/0.087/0.266,
      xypp/-0.033/0.452, xypm/-0.434/0.081,
      xymp/0.583/0.120, xymm/0.151/-0.208,
      yzpp/0.251/-0.019, yzpm/0.312/0.559,
      yzmp/-0.180/-0.355, yzmm/-0.093/0.201,
      zxpp/-0.295/-0.037, zxpm/0.331/-0.344,
      zxmp/-0.194/0.531, zxmm/0.383/0.221}{
      \coordinate (\vertex) at (\u,\v);
    }

    \draw[wire] (outerxm) -- (outerym) -- (outerzm) -- cycle;
    \foreach \s in {p,m}{
      \draw[rear wire]
        (xym\s) -- (outerxm) -- (zx\s m)
        (yzm\s) -- (outerym) -- (xy\s m)
        (zxm\s) -- (outerzm) -- (yz\s m);
    }

    \foreach \a in {p,m} \foreach \b in {p,m}{
      \foreach \c in {p,m}{
        \draw[rear wire] (xy\a\b) -- (yz\b\c) -- (zx\c\a) -- cycle;
      }
      \draw[rear wire]
        (innerx\a) -- (xy\a\b) -- (innery\b)
        (innery\a) -- (yz\a\b) -- (innerz\b)
        (innerz\a) -- (zx\a\b) -- (innerx\b);
      \draw[wire]
        (innerx\a) -- (innery\b)
        (innery\a) -- (innerz\b)
        (innerz\a) -- (innerx\b);
    }

    \draw[black,line width=.49pt]
      (outerxp) -- (outerzm) -- (outeryp) -- (outerxm)
      -- (outerzp) -- (outerym) -- cycle;
    \foreach \s in {p,m}{
      \draw[near wire]
        (xyp\s) -- (outerxp) -- (zx\s p)
        (yzp\s) -- (outeryp) -- (xy\s p)
        (zxp\s) -- (outerzp) -- (yz\s p);
    }

    \foreach \a in {p,m}{
      \foreach \axis in {x,y,z}
        \fill (inner\axis\a) circle[radius=.012];
      \foreach \b in {p,m} \foreach \pair in {xy,yz,zx}
        \fill (\pair\a\b) circle[radius=.017];
    }
    \draw[black,line width=.78pt,preaction={draw=white,line width=1.30pt}]
      (outerxp) -- (outeryp) -- (outerzp) -- cycle;
    \foreach \axis in {x,y,z}{
      \fill (outer\axis m) circle[radius=.026];
      \fill (outer\axis p) circle[radius=.038];
    }
  \end{scope}

  \node[font=\Large] at (10.30,0) {?};
  \foreach \pos in {1.85,5.55,9.15}{
    \node[font=\normalsize] at (\pos,0) {$\to$};
  }
  \foreach \pos/\dimension in {0/{n=2},3.70/{n=3},7.40/{n=4},10.30/{n\to\infty}}{
    \node[anchor=base] at (\pos,-1.57) {$\dimension$};
  }
  \node[anchor=north] at (0,-1.70) {\strut hexagon};
  \node[anchor=north] at (3.70,-1.70) {\strut truncated octahedron};
  \node[anchor=north] at (7.40,-1.70) {\strut $24$-cell};
\end{tikzpicture}
\caption{The three known minimizers in dimensions 2, 3 and 4 are all lattice Voronoi cells.
Voronoi conjectured that every parallelohedron is an affine image of such a cell \cite{Voronoi1908,Voronoi1909}.
A recent breakthrough came in dimension five, where Garber proved the conjecture \cite{Garber2025}.}
\label{fig:optimal-shapes}
\end{figure}

In the plane, Hales proved the honeycomb conjecture \cite{Hales2001}.
Fejes T\'oth had earlier proved the conjecture with the convexity hypothesis \cite{FejesToth1943}.
As a corollary, the regular hexagon minimizes perimeter among planar parallelohedra of fixed area, since every such parallelohedron tiles the plane by translations.

In three dimensions, Bezdek's conjecture that the Archimedean truncated octahedron minimizes surface area at fixed volume among parallelohedra \cite{Bezdek2006} was proved by Hales and Song \cite{HalesSong2026} and independently by Cesaroni and Novaga \cite{CesaroniNovagaGlobal2026}.
The local theory was developed by Cesaroni and Novaga for lattice Voronoi cells \cite{CesaroniNovagaLocal2026} and by Song for parallelohedra \cite{Song2026}.
L\'angi had earlier proved that the Archimedean truncated octahedron also minimizes mean width for fixed volume \cite{Langi2022}.
Cesaroni and Novaga proved that the regular rhombic dodecahedron minimizes surface area among parallelohedra of fixed inradius \cite{CesaroniNovagaRhombic2026}.

In four dimensions, Novaga and Song proved that the regular $24$-cell minimizes surface area among parallelohedra of fixed volume or fixed inradius \cite{NovagaSong2026}, while 
the problem seems currently out of reach in dimensions five or higher.

Could the optimal isoperimetric quotient follow a simple law as the dimension tends to infinity, even when the minimizing shapes remain unknown?

We answer this question affirmatively by proving that at equal volume, the smallest possible surface area of a parallelohedron approaches $\sqrt{2}$ times that of the ball.
The same method produces convex $\Z^n$-tiles all of whose hyperplane projections are small, which answers a question of Naor and Regev \cite{NaorRegev2023}.

For a convex body $K\subset\R^n$, write $\vol(K)$ for its volume and $\area(\partial K)$ for its surface area. Note its surface area is the $(n-1)$-dimensional Hausdorff measure of its boundary.
Its isoperimetric quotient is
\[
 \iq(K)=\frac{\area(\partial K)}{\vol(K)^{(n-1)/n}}.
\]
Let $B_1^n$ denote the Euclidean unit ball in $\R^n$.
Let $\omega_n=\vol(B_1^n)$ and $b_n=n\omega_n^{1/n}$.
Thus $b_n$ is the isoperimetric quotient of a ball.
For a full rank lattice $\Lambda\subset\R^n$, let $\Vor\Lambda$ denote its ordinary Euclidean Voronoi cell centered at the origin.
We use the infima
\begin{align*}
 \PP_n&=\inf\{\iq(K):K\subset\R^n\text{ is a parallelohedron}\},\\
 \AV_n&=\inf\{\iq(K):K\text{ is an affine image of a lattice Voronoi cell}\},\\
 \VV_n&=\inf\{\iq(\Vor\Lambda):\Lambda\subset\R^n\text{ is a full rank lattice}\}.
\end{align*}
The Euclidean isoperimetric inequality and the inclusions among these classes give
\begin{equation}\label{eq:hierarchy}
 b_n\le\PP_n\le\AV_n\le\VV_n.
\end{equation}

\begin{theorem}\label{thm:main}
There is an absolute constant $C>0$ such that every $n$-dimensional parallelohedron $K$, $n\ge2$, satisfies
\begin{equation}\label{eq:lower-rate}
 \frac{\iq(K)}{b_n}\ge\sqrt2\,(1-Cn^{-1/2}).
\end{equation}
Moreover, for each $Q_n\in\{\PP_n,\AV_n,\VV_n\}$,
\begin{equation}\label{eq:limits}
 \frac{Q_n}{b_n}=\sqrt2+O(n^{-1/2}).
\end{equation}
In particular, $\PP_n/\VV_n=1+O(n^{-1/2})$.
\end{theorem}

\subsection{Ideas of the proof}

Both bounds use the same fact about two balls.
If their radii are $r$ and their centers are at distance $d\le 2r$, the parallelogram identity places their intersection inside the ball of radius $\sqrt{r^2-d^2/4}$ centered halfway between them.

For the lower bound, center a Gaussian at a point $y$ and record the mass it assigns to each tile.
We first keep finitely many tiles and group the rest into one region.
After rescaling to a standard Gaussian, the multi-bubble theorem of Milman and Neeman \cite{MN} bounds the Gaussian perimeter of this partition from below by that of a simplicial cluster with the same masses.
The rule defining this cluster is simple: add fixed offsets to the Gaussian coordinates and select the largest adjusted coordinate. Choose the offsets so that label $i$ has probability $v_i$.
The perimeter is $1/\sqrt2$ times the expected value of the selected coordinate before its offset is added.
The factor $1/\sqrt2$ is the absolute value of the $i$-th component of the unit normal $(e_j-e_i)/\sqrt2$ pointing from cell $i$ to cell $j$.

Among all rules with the same label probabilities, the simplicial rule maximizes the expected value of the selected coordinate. To bound this expectation from below one tile at a time, we turn Gaussian coordinate $i$ into an exponential waiting time of mean $1/v_i$.
Averaging over Gaussian centers in one lattice period turns this sum into an integral over $\R^n$ involving a single tile.
Because the contribution is concave in the tile mass, rearrangement shows that replacing the tile by a ball of equal volume can only decrease this integral.

The smoothed density $q(y)$ is the probability that a Gaussian displacement from \(y\) lands in the ball.
We compare the clock contribution at mass $q(y)$ with the Gaussian perimeter of a halfspace of that mass.
Radial integration estimates the main term, while the expected overlap of two independently displaced copies of the ball controls the loss in this comparison.
At a single smoothing scale, the two-ball estimate gives both the uniform error in \eqref{eq:lower-rate} and the explicit lower bound in Theorem~\ref{thm:finite}.
Proposition~\ref{prop:race-optimum} identifies the optimal scale and the first correction to the smoothed ball bound.

For the upper bound, choose a point uniformly in a lattice Voronoi cell of volume one and follow the ray from the origin through that point until it hits a facet.
The probability of hitting a given facet equals the volume of the pyramid it spans with the origin. This volume is $Ah/n$ if the facet has area $A$ and its plane is at distance $h$ from the origin.
It follows that the surface area of the cell is $n$ times the mean of $1/h$ for the facet hit by the ray.
Averaging over lattices with Siegel mean value theorem \cite{Siegel1945} bounds the contribution of facets near the origin by the volume of a spherical cap.
The cap is half the overlap of two equal balls, so the same two-ball estimate gives the matching upper bound in Proposition~\ref{prop:mean-upper}.
Corollary~\ref{cor:typical-lattice} quantifies how closely a typical lattice attains it.

For the same random facet, the second moment $\E(1/h^2)$ controls every hyperplane projection by \eqref{eq:projection-height}.
We choose a lattice $B\Z^n$ whose Voronoi cell has all projection volumes $O(1)$ and whose basis matrix satisfies $\det B=1$ and $\lVert B\rVert_{\mathrm{op}}=O(\sqrt{\log n})$.
Applying $B^{-1}$ gives a convex $\Z^n$-tile whose largest projection volume is $O(\sqrt{\log n})$, as asserted in Theorem~\ref{thm:integer-projections}.

\subsection{Related work}
The high-dimensional isoperimetric problem for translational tiles depends crucially on
whether convexity is imposed and whether the lattice of translations is prescribed.
Without the convexity constraint, Kindler, O'Donnell, Rao and Wigderson constructed
fundamental domains for $\R^n/\Z^n$ with surface area $O(\sqrt n)$
\cite{KindlerEtAl2008}.
Their construction is not convex; the same paper also discusses small area lattice
Voronoi cells, by a different argument.
Alon and Klartag gave a short proof of an $O(\sqrt n)$ bound for a related
problem concerning sets that intersect every noncontractible closed curve in the
flat torus $\R^n/\Z^n$ \cite{AlonKlartag2009}.
In another direction, Braverman and Minzer showed that requiring a $\Z^n$-tile
to be invariant under coordinate permutations entails a substantially larger
minimum surface area, of order $n/\sqrt{\log n}$, even without assuming convexity
\cite{BravermanMinzer2021}.

For convex translational tiles, allowing the lattice to vary makes the order
$\sqrt n$ accessible through classical results in the geometry of numbers.
Rogers' lattice packing estimates give volume-one lattice Voronoi cells with
surface area $O(\sqrt n)$ \cite{Rogers1950}.
More precisely, Butler's simultaneous packing and covering theorem yields
lattice Voronoi cells with isoperimetric quotient at most
$(2+o(1))b_n$, as noted in \cite{KindlerEtAl2008,Butler1972}.
A related upper bound, based on lattice packing estimates, was obtained by
Cesaroni and Novaga in the setting of periodic perimeter-minimizing foams
\cite{CesaroniNovaga2024}.
These results determined the correct order of growth but not the sharp
asymptotic coefficient. Theorem~\ref{thm:main} identifies this coefficient
as $\sqrt2$ and proves a matching lower bound for \emph{every} parallelohedron,
not only for suitably chosen lattice Voronoi cells.

The problem becomes more restrictive when the lattice is fixed.
Naor and Regev constructed convex $\Z^n$-tiles with surface area $n^{1/2+o(1)}$.
Question~2 in~\cite{NaorRegev2023} asks whether such tiles can have surface area $O(\sqrt n)$.
Question~14 in the same paper asks whether there are convex $\Z^n$-tiles whose orthogonal projections onto \emph{every} hyperplane have $(n-1)$-dimensional volume $n^{o(1)}$.
Theorem~\ref{thm:integer-projections} answers Question~14 by giving an $O(\sqrt{\log n})$ bound on the largest projection volume.
Cauchy's projection formula then gives a surface area bound of $O(\sqrt{n\log n})$ for the same tiles.
Ball, B\"or\"oczky and Naor proposed a polynomial-time optimization oracle with a constant approximation factor as an efficient description of such convex tiles in high dimensions \cite{BallBoroczkyNaor2026}.
Corollary~\ref{cor:integer-oracle} supplies an oracle for the tiles constructed here, with an approximation factor tending to $1/4$.

There is also some literature on isoperimetric problems for nonconvex tiles.
Choe studied existence and regularity of fundamental domains of least area in covering spaces \cite{Choe1989}.
Cesaroni and Novaga established existence of perimeter-minimizing periodic
tilings for a fixed lattice \cite{CesaroniNovaga2024Partitions} and when the
lattice varies at fixed cell volume \cite{CesaroniNovaga2024};
they also investigated the corresponding problem with a prescribed inradius
\cite{CesaroniNovagaInradius2025}.
These problems permit nonconvex cells and therefore differ from the
minimization over parallelohedra considered here.

Related results in fixed dimension include the optimality of regular hexagons
for several functionals among planar convex translational tiles
\cite{CesaroniFragalaNovaga2026} and the honeycomb theorem of Lángi and Wang for normed planes \cite{LangiWang2026}.
Quantitative stability near the Euclidean isoperimetric optimum was proved by
Fusco, Maggi and Pratelli \cite{FuscoMaggiPratelli2008}. Related stability results for planar honeycombs and finite clusters appear in
\cite{CarocciaMaggi2016,CarocciaDeMasonMaggi2026}.
These concern a different regime from the asymptotic factor $\sqrt2>1$ in
Theorem~\ref{thm:main}.
Finally, Ball, B\"or\"oczky and Naor obtained sharp-order isoperimetric
bounds for polytopes with a prescribed number of facets
\cite{BallBoroczkyNaor2026}, while Jo\'os and L\'angi solved related
minimization problems for zonotopes generated by $n$ or $n+1$ segments
\cite{JoosLangi2023}.

\section{The asymptotic isoperimetric quotient}\label{sec:asymptotic}

\subsection{A Gaussian comparison for the tiling}\label{sec:clusters}
Write $\gamma_N$ for standard Gaussian measure on $\R^N$, and $\varphi,\Phi,\bar\Phi=1-\Phi$ for the one-dimensional density, distribution function and upper tail.
The Gaussian isoperimetric profile is $\J(u)=\varphi(\Phi^{-1}(u))$ for $0<u<1$, with $\J(0)=\J(1)=0$.

We seek a lower bound on the Gaussian area of all the interfaces that is a sum of contributions from the individual tile masses.
Averaging this bound over the Gaussian center will give a lower bound for the surface area of one tile.

A finite Gaussian cluster is a partition, up to null sets, into Borel sets $\Omega_i$ of finite Gaussian perimeter.
Its perimeter $\mathcal P_{\gamma_N}(\Omega)=\frac12\sum_iP_{\gamma_N}(\Omega_i)$ counts each interface once.
On a polyhedron, $P_{\gamma_N}$ is ordinary surface area weighted by the Gaussian density.
For a positive probability vector $v=(v_1,\dots,v_q)$, let $I_q(v)$ be the perimeter of the translated simplicial cluster of masses $v$ in $E=\{x\in\R^q:\sum_i x_i=0\}$.
Its cells are obtained by translating the partition of $E$ according to which coordinate is largest.
The Gaussian multi-bubble theorem of Milman and Neeman gives
\begin{equation}\label{eq:MN}
 \mathcal P_{\gamma_N}(\Omega)\ge I_q(v),
\end{equation}
where $v_i=\gamma_N(\Omega_i)$.
The theorem applies directly for $2\le q\le N+1$.
For larger finite $q$, take the product of the cluster with a Gaussian space of dimension $q-1-N$.
This preserves its masses and perimeter \cite{MN}.

\begin{lemma}\label{lem:scores}
Let $\xi_1,\dots,\xi_q$ be independent standard normals.
There are offsets $b_i$ for which $I=\argmax_i(\xi_i+b_i)$ has law $v$, and
\begin{equation}\label{eq:model-score}
 \E\xi_I=\sqrt2\,I_q(v).
\end{equation}
For every choice of label $J$ with the same probabilities,
\begin{equation}\label{eq:variational}
 \E\xi_J\le\E\xi_I .
\end{equation}
\end{lemma}
\begin{proof}
Choose the translate of the simplicial cluster with masses $v_i$ \cite{MN}.
Orthogonal projection onto $E$ subtracts the coordinate mean.
It preserves all coordinate differences, so the inverse images of the cells have the form $C_i=\{x:x_i+b_i\ge x_j+b_j\text{ for all }j\}$.
Thus $\Prb(I=i)=v_i$.
Ties have probability zero.

On an interface between $C_i$ and $C_j$ the outward normal from $C_i$ is $(e_j-e_i)/\sqrt2$.
Gaussian integration by parts thus gives $\int_{C_i}x_i\,d\gamma_q=P_{\gamma_q}(C_i)/\sqrt2$.
Summing over $i$ proves \eqref{eq:model-score}.
Finally $\xi_J+b_J\le\xi_I+b_I$ pointwise and $\E b_J=\E b_I$, proving \eqref{eq:variational}.
\end{proof}

For prescribed masses $v_1,\dots,v_q$, let clock $i$ ring at time $T_i=-\log\Phi(\xi_i)/v_i$.
These waiting times are independent exponentials with rates $v_i$, and we let $J$ denote the index of the first clock to ring.
Given $\xi_i=x$, the probability that all the other clocks ring later is $\Phi(x)^{(1-v_i)/v_i}$, since their rates sum to $1-v_i$.
Integrating against $\varphi(x)$ gives
\[
 \Prb(J=i,\,\xi_i\le t)=\int_{-\infty}^{t}\varphi(x)\Phi(x)^{1/v_i-1}\,dx
 =v_i\Phi(t)^{1/v_i}.
\]
Thus label $i$ has probability $v_i$, and its selected coordinate has conditional distribution function $\Phi^{1/v_i}$.

For $z\ge0$, let $m(z)$ be the mean of the distribution with cumulative distribution function $\Phi^{e^z}$.
When $e^z$ is an integer this is the mean of the largest of $e^z$ independent standard normals.
A label of mass $u$ therefore contributes $F(u)=u\,m(-\log u)$ to the expected score, and we define $F(0)=0$.
This contribution is nonnegative because $\Phi^{1/u}\le\Phi$ for $0<u\le1$.
The score comparison \eqref{eq:variational} now gives the desired sum,
\begin{equation}\label{eq:race-scores}
 \sqrt2\,I_q(v)\ge\E\xi_J=\sum_{i=1}^q F(v_i).
\end{equation}
For equal masses the first clock selects the largest Gaussian coordinate, and equality holds.
Scale the parallelohedron $K$ to volume one and choose a lattice tiling $\{K+\ell:\ell\in\Lambda\}$.
For $\sigma>0$ write
\begin{align*}
 \gamma_\sigma(x) &=(2\pi\sigma^2)^{-n/2}e^{-\lVert x\rVert^2/(2\sigma^2)}, \\
 p &=\ind_K*\gamma_\sigma, \\
 p_\ell(y) &=p(y-\ell).
\end{align*}
The numbers $p_\ell(y)$ are the tile masses under $N(y,\sigma^2 I_n)$ and sum to one.
The total interface area seen by this Gaussian is
\[
 A(y)=\frac12\sum_{\ell\in\Lambda}
       \int_{\partial(K+\ell)}\gamma_\sigma(x-y)\,d\mathcal H^{n-1}(x).
\]
It is finite by boundedness of $K$ and Gaussian decay.

\begin{proposition}\label{prop:window}
For every $y$ and $\sigma>0$,
\begin{equation}\label{eq:window}
 \sigma A(y)\ge\frac1{\sqrt2}\sum_{\ell\in\Lambda}F(p_\ell(y)).
\end{equation}
Consequently, if $X$ is uniform on $K$, $Z$ is an independent standard Gaussian vector and $L=-\log p(X+\sigma Z)$, then
\begin{equation}\label{eq:race-tile}
 \area(\partial K)\ge\frac{\sqrt2}{\sigma}\int_{\R^n}F(p(y))\,dy
 =\frac{\sqrt2}{\sigma}\E m(L).
\end{equation}
\end{proposition}
\begin{proof}
Keep finitely many tiles and merge all others into their complement.
The resulting cluster has finite perimeter and all its interfaces belong to the original tiling.
Indeed, the retained union is a finite union of polytopes and grouping cells only removes pieces of their boundaries up to sets of surface measure zero \cite{Maggi2012}.
The map $x\mapsto(x-y)/\sigma$ preserves cluster masses and multiplies weighted perimeter by $\sigma$.
By \eqref{eq:MN} and \eqref{eq:race-scores}, $\sqrt2\,\sigma A(y)$ bounds the sum of $F(p_\ell(y))$ over the retained tiles, since the contribution of the merged cell is nonnegative.
Letting the finite sets increase to $\Lambda$ proves \eqref{eq:window}.

Let $D$ be a fundamental domain of $\Lambda$.
Its translates cover $\R^n$ and $\vol(D)=1$.
Integrating over $y\in D$ and translating the summands gives
\[
 \int_D A(y)\,dy=\frac12\area(\partial K)
 \quad\text{and}\quad
 \int_D\sum_{\ell\in\Lambda}F(p_\ell(y))\,dy
       =\int_{\R^n}F(p(y))\,dy.
\]
The first identity uses $\int\gamma_\sigma=1$.
Since $Y=X+\sigma Z$ has density $p$ and $F(p)=p\,m(-\log p)$, the second gives \eqref{eq:race-tile}.
\end{proof}

Applying Gaussian isoperimetry to each tile separately counts every interface twice and gives only $2\sigma A(y)\ge\sum_\ell\J(p_\ell(y))$ \cite{Borell1975}.

\subsection{The uniform lower bound for parallelohedra}\label{sec:ball}

The Gaussian cluster profile $I_q$ is strictly concave in the cell masses \cite{MN}, so Jensen’s inequality bounds its average over Gaussian centers from above. For our lower bound, the clocks give a sum of individual tile contributions, and the concavity of $F$ lets us compare the resulting integral with that of a ball.
Integration of the two tails gives
\[
 F(u)=u\int_0^\infty
       \bigl(1-\Phi(t)^{1/u}-\bar\Phi(t)^{1/u}\bigr)\,dt.
\]
The integrand is nonnegative since $1/u\ge1$.
For $d>0$ the function $u\mapsto u e^{-d/u}$ is convex, so $F$ is concave.

\begin{proposition}\label{prop:race-ball}
Let $B$ be the centered ball of volume one and radius $r_n=\omega_n^{-1/n}$, and define $q=\ind_B*\gamma_\sigma$.
Every parallelohedron $K\subset\R^n$ satisfies
\begin{equation}\label{eq:race-comparison}
 \frac{\iq(K)}{\sqrt2\,b_n}\ge D_n(v):=
 \frac1{\sigma b_n}\int_{\R^n}F(q(y))\,dy
 \quad\text{where}\quad
 v=\frac{n\sigma^2}{r_n^2}>0 .
\end{equation}
\end{proposition}
\begin{proof}
Normalize $K$ to volume one and let $p=\ind_K*\gamma_\sigma$.
For every measurable set $E$ of finite volume, let $E^*$ be the centered ball of the same volume.
Since $\gamma_\sigma$ is radial and decreasing, the Riesz rearrangement inequality gives $\int_E p\le\int_{E^*}q$.
Write $u_+=\max\{u,0\}$ for the positive part.
Taking $E=\{p>t\}$ for $0<t<1$ and subtracting $t\vol(E)$ from both sides gives
\[
 \int_{\R^n}(p-t)_+\le\int_{E^*}(q-t)\le\int_{\R^n}(q-t)_+ .
\]
Both densities have integral one, so $\int\min(p,t)\ge\int\min(q,t)$.
Every concave polygonal function on $[0,1]$ that vanishes at zero is a linear function plus a sum of terms $c_t\min(u,t)$, where $c_t\ge0$ is the drop in slope at the vertex $t$.
By concavity, the polygonal interpolants of $F$ on successive dyadic subdivisions increase to $F$ at every interior point, and they agree with $F$ at both endpoints.
Monotone convergence therefore gives $\int F(p)\ge\int F(q)$, and \eqref{eq:race-tile} gives the claim.
Wang and Madiman give the general convolution comparison as Theorem~VII.4 in \cite{WangMadiman2014}.
\end{proof}

For a Gaussian halfspace of mass $u$, the mean coordinate in its inward direction is $\J(u)/u$.
The next estimate compares this mean with the score of a label of the same mass.

\begin{lemma}\label{lem:score-gap}
For $0<u<1$,
\[
 0<\frac{\J(u)-F(u)}{u}\le\frac{2}{\sqrt{3/2-\log u}}.
\]
\end{lemma}
\begin{proof}
Write $u=\bar\Phi(c)$, and let $M$ have distribution function $\Phi^{1/u}$, so $\E M=F(u)/u$.
The density of $M$ multiplied by $u$ is $\varphi\Phi^{1/u-1}$.
It has mass $u$ and is bounded by $\varphi$.
Among such densities, the upper normal tail of mass $u$ uniquely maximizes the first moment, giving $F(u)<\J(u)$.

Let $\lambda=\varphi(c)/u$ be the mean of that tail and write $\Delta=\lambda-\E M$.
Since $\Phi(M)^{1/u}$ is uniform on $(0,1)$, concavity of $\log\Phi$ gives $\Phi(\E M)\ge e^{-u}\ge1-u$, hence $\E M\ge c$.
For $c\le0$, the claim follows from $u\ge1/2$, $\E M\ge0$ and $\lambda\le2\varphi(0)<1$.
For $c>0$, let $Z$ be standard normal.
Integration by parts gives $\E(Z^2\mid Z\ge c)=1+c\lambda$, so nonnegativity of the conditional variance gives $\lambda(\lambda-c)\le1$.
Thus the mean excess $e=\lambda-c$ satisfies $e(c+e)\le1$, and
\[
 \Delta\le e\le\frac{2}{c+\sqrt{c^2+4}}\le\frac{2}{\sqrt{c^2+4}}.
\]
Integrating $\varphi$ from $c$ to $c+1$ gives $u\ge\varphi(c+1)$, hence $-\log u\le c^2+2$, which proves the claim.
\end{proof}

\begin{proof}[Proof of the uniform lower bound \eqref{eq:lower-rate}]
Let $Y_B$ have density $q$ and write $a=nv=(\sigma b_n)^2$.
We first estimate the integral in $D_n(v)$ with the halfspace profile $\J(q)$ in place of $F(q)$.

Differentiating the Gaussian kernel gives $\sigma\nabla q(y)=\E(Z\ind_{\{y+\sigma Z\in B\}})$ for a standard Gaussian vector $Z$.
The event has probability $q(y)$, and a halfspace of this mass maximizes the first moment in any fixed direction, so $\sigma\lVert\nabla q(y)\rVert\le\J(q(y))$.
Since $q$ is radial and decreasing, integration by parts in the radius gives $\int\lVert\nabla q\rVert=(n-1)\E(1/\lVert Y_B\rVert)$.
Now $Y_B$ is the sum of a uniform point of $B$ and an independent Gaussian of covariance $\sigma^2 I_n$, so $\E\lVert Y_B\rVert^2=r_n^2(n/(n+2)+v)$.
Jensen's inequality and $b_n=n/r_n$ give
\begin{equation}\label{eq:ball-profile}
 \frac{1}{\sigma b_n}\int_{\R^n}\J(q(y))\,dy
 \ge\frac{n-1}{b_n\sqrt{\E\lVert Y_B\rVert^2}}
 =\frac{1-1/n}{\sqrt{n/(n+2)+v}}.
\end{equation}

To bound the score loss, observe that $\E q(Y_B)=\int q^2$ is the expected volume of overlap of two independently displaced copies of $B$.
Their relative displacement has the law of $h=\sqrt2\,\sigma Z$.
For $x\in B\cap(B+h)$, the parallelogram identity gives
\[
 \left\lVert x-\frac h2\right\rVert^2\le r_n^2-\frac{\lVert h\rVert^2}{4},
\]
which bounds the overlap by the volume of a smaller ball.
Since $B$ has volume one, this gives
\[
 \E q(Y_B)=\E\vol\bigl(B\cap(B+h)\bigr)
 \le \E\exp\!\left(-\frac{n\lVert h\rVert^2}{8r_n^2}\right)
 =(1+v/2)^{-n/2}.
\]
The function $g(u)=(3/2-\log u)^{-1/2}$ is increasing and concave on $(0,1]$, since its second derivative has the sign of $\log u$.
Lemma~\ref{lem:score-gap} and Jensen's inequality thus give
\begin{equation}\label{eq:ball-score-loss}
\begin{aligned}
 \frac1{\sigma b_n}\int_{\R^n}\bigl(\J(q(y))-F(q(y))\bigr)\,dy 
 &\le \frac{2}{\sqrt a}\,g\bigl(\E q(Y_B)\bigr) \\
 &\le \frac{2}{\sqrt{a\left(\dfrac{3}{2}+\dfrac{n}{2}\log(1+\dfrac{v}{2})\right)}}.
\end{aligned}
\end{equation}
At $v=\sqrt{8/n}$, the right side of \eqref{eq:ball-profile} is $1-O(n^{-1/2})$ and the loss in \eqref{eq:ball-score-loss} is $O(n^{-1/2})$.
Subtracting and using \eqref{eq:race-comparison} proves \eqref{eq:lower-rate}.
\end{proof}

\subsection{The matching upper bound and typical Voronoi cells}\label{sec:upper}

Let $L$ be a covolume-one lattice chosen according to Haar probability $\mu_n$ on the space $\mathcal L_n$ of such lattices.
Write $K=\Vor L$ and $R_n=\iq(K)/b_n$.
Each facet of $K$ lies on the perpendicular bisector of $0$ and a lattice point $v\ne0$.
Its height above the origin is $\lVert v\rVert/2$.
The question is how much surface area belongs to facets of small height.

Choose $X$ uniformly in $K$ and extend the ray from $0$ through $X$ to its first boundary point $Y$.
Almost surely $Y$ lies on a unique facet.
Denote its height by $r_n H$ and let $S=(\lVert Y\rVert/r_n)^n$, the volume of the ball with radius $\lVert Y\rVert$.
A facet element $dA$ and the origin span a pyramid of volume $r_n H\,dA/n$.
This is the probability that the ray ends in $dA$.
In particular, since $b_n=n/r_n$,
\begin{equation}\label{eq:cone-identity}
 \E(H^{-1}\mid L)=R_n .
\end{equation}

We use Siegel's mean value theorem \cite{Siegel1945} in two ways.
For every nonnegative integrable $f$ it states that
\begin{equation}\label{eq:siegel}
 \E\sum_{v\in L\setminus\{0\}} f(v)=\int_{\R^n}f(x)\,dx .
\end{equation}
The first application bounds the chance that the ray ends far away.
For a centered ball $J$ of volume $V$, let $N$ count the points of a uniform random coset of $L$ in $J$.
Averaging over the coset gives $\E(N\mid L)=V$.
Counting ordered pairs with difference $v\ne0$ gives
\[
 \E(N(N-1)\mid L)
 =\sum_{v\in L\setminus\{0\}}\vol\bigl(J\cap(J+v)\bigr).
\]
By \eqref{eq:siegel} the latter sum has mean $V^2$, so $\Prb(N=0)\le\Var(N)/V^2=1/V$.
The point $X$ represents a uniform coset by its point of least norm.
Thus $\Prb(X\notin J)=\Prb(N=0)$.
By radial integration, conditional on $Y$, the ratio $(\lVert X\rVert/\lVert Y\rVert)^n$ is uniform on $[0,1]$.
Thus $\Prb(X\notin J\mid Y)=(1-V/S)_+$.
Taking $V=M/2$ shows that
\begin{equation}\label{eq:cone-far}
 \Prb(S>M)\le2\Prb(X\notin J)\le4/M.
\end{equation}

The second application counts the possible facets within the ball of volume $M$.
Write $\rho=M^{1/n}$ and $g_n=\omega_{n-1}/\omega_n$.
A facet of normalized height $h$ cuts from the ball of radius $r_n\rho$ a section of area at most $\omega_{n-1}r_n^{n-1}(\rho^2-h^2)_+^{(n-1)/2}$.
Multiply by the pyramid weight $r_n h/n$ and sum over all possible labels $v$ of length $2r_n h$.
For every nonnegative measurable function $f$, Siegel's formula in polar coordinates gives
\begin{equation}\label{eq:facet-measure}
 \E\bigl(f(H)\ind_{\{S\le M\}}\bigr)
 \le 2g_n\int_0^\rho f(h)\,h
       \bigl(2h\sqrt{\rho^2-h^2}\bigr)^{n-1}\,dh .
\end{equation}
The factor $h\sqrt{\rho^2-h^2}$ is largest at $h=\rho/\sqrt2$.
This is where the constant in the upper bound enters.

\begin{proposition}\label{prop:mean-upper}
For an absolute constant $C$ and every $n\ge2$,
\begin{equation}\label{eq:upper-rate}
 \frac{\VV_n}{b_n}\le\E R_n\le\sqrt2+\frac C{\sqrt n}.
\end{equation}
\end{proposition}
\begin{proof}
Define $G(t)=\Prb(H\le t)$.
For $0\le t\le\rho$, substitute $z=2h^2/\rho^2-1$ in \eqref{eq:facet-measure} with $f=\ind_{[0,t]}$.
If $U$ is uniform in the unit ball, its first coordinate has density $g_n(1-z^2)^{(n-1)/2}$ on $[-1,1]$.
Hence
\[
 \Prb(H\le t \text{ and } S\le M)
 \le \frac12 M^2\Prb\bigl(U_1\le 2t^2/\rho^2-1\bigr).
\]
For $0\le d\le1$, reflection in $x_1=d$ splits $B_1^n\cap(B_1^n+2de_1)$ into two congruent caps, one of which is $\{x\in B_1^n:x_1\ge d\}$.
The two-ball inclusion from Section~\ref{sec:ball} places the intersection in a ball of radius $\sqrt{1-d^2}$, so $\Prb(U_1\le-d)=\Prb(U_1\ge d)\le(1-d^2)^{n/2}/2$.
For $0\le t\le1/\sqrt2$, let $d=1-2t^2$.
Whenever $M\ge1$, \eqref{eq:cone-far} and the cap bound give
\[
 G(t)\le 4/M+(M^2/4)e^{-nd^2/2}.
\]
Choosing $M=e^{nd^2/6}$ gives the first of the two estimates
\begin{equation}\label{eq:height-tails}
 G(t)\le\frac{17}{4}e^{-n(1-2t^2)^2/6}
 \quad\text{and}\quad
 G(t)\le\frac{(2t)^n}{2}
\end{equation}
for $0\le t\le1/\sqrt2$.
The second follows from \eqref{eq:siegel} because $H\le t$ requires two opposite lattice vectors in the ball of radius $2r_n t$.

Since $1-2/u^2\ge (u-\sqrt2)/3$ for $\sqrt2\le u\le3$, these two bounds give
\[
 \Prb(H^{-1}>u)\le
 \begin{cases}
  (17/4)e^{-n(u-\sqrt2)^2/54} & \text{for } \sqrt2\le u\le3,\\
  (2/u)^n/2 & \text{for } u>3 .
 \end{cases}
\]
A conditional mean can only decrease the expected positive part.
Thus \eqref{eq:cone-identity} and integration from $a$ to infinity give
\begin{equation}\label{eq:mean-excess}
 \E(R_n-a)_+\le\E(H^{-1}-a)_+
 \le\frac C{\sqrt n}e^{-cn(a-\sqrt2)^2}
 \quad\text{for}\quad
 \sqrt2\le a\le\sqrt2+\frac12 .
\end{equation}
Here the Gaussian integral uses $(u-\sqrt2)^2\ge(u-a)^2+(a-\sqrt2)^2$ for $u\ge a$.
The integral above $3$ is $3(2/3)^n/\bigl(2(n-1)\bigr)$.
Take $a=\sqrt2$ and use $\VV_n/b_n\le\E R_n$.
\end{proof}

\begin{proof}[Conclusion of the proof of Theorem~\ref{thm:main}]
Every normalized quotient is at least $\sqrt2-O(n^{-1/2})$ by \eqref{eq:lower-rate}.
Equation \eqref{eq:upper-rate} gives the opposite bound for the smallest ordinary Voronoi quotient.
The inclusions \eqref{eq:hierarchy} place all three infima between these bounds.
\end{proof}

\begin{corollary}\label{cor:typical-lattice}
For a Haar-random lattice $L$ of covolume one, let $R_n=\iq(\Vor L)/b_n$.
There are absolute constants $C,c>0$ such that
\[
 \E|R_n-\sqrt2|\le\frac C{\sqrt n}
 \quad\text{and}\quad
 \Prb\!\left(|R_n-\sqrt2|>\frac{s}{\sqrt n}\right)
       \le C e^{-cs^2}
\]
for $1\le s\le\sqrt n$.
\end{corollary}
\begin{proof}
The uniform lower bound gives $(\sqrt2-R_n)_+\le D/\sqrt n$ for an absolute $D$.
Add this to \eqref{eq:mean-excess} at $a=\sqrt2$ for the first estimate.
For the second take $a=\sqrt2+s/(2\sqrt n)$.
An upper deviation of $s/\sqrt n$ requires $(R_n-a)_+>s/(2\sqrt n)$, so Markov's inequality and \eqref{eq:mean-excess} bound its probability by $(2C/s)e^{-cs^2/4}$.
A lower deviation is possible only for $s<D$, where an increase of the constant suffices.
\end{proof}

For the typical Poisson--Voronoi cell $Z_n$, H\"orrmann, Hug, Reitzner and Th\"ale obtained the same limit $\sqrt2$ for the ratio $\E\,\area(\partial Z_n)/\bigl(b_n(\E\,\vol(Z_n))^{(n-1)/n}\bigr)$, see Theorem~3.20 in~\cite{HorrmannHugReitznerThale2015}.
In the lattice model each cell has volume one.

\subsection{The sharp correction to the smoothed ball bound}\label{sec:sharp-ball}
The proof of \eqref{eq:lower-rate} bounded the two losses from smoothing by $O(v)$ and $O(1/(nv))$ when $1/n\le v\le1$.
We now determine their leading coefficients.
Return to the density $q=\ind_B*\gamma_\sigma$ in \eqref{eq:race-comparison}, let $Y_B$ have density $q$, and write $L=-\log q(Y_B)$ and $a=nv=(\sigma b_n)^2$.
Then $D_n(v)=\E m(L)/\sqrt a$.
With $H=(n/2)\log(1+v)$, the mean and variance of $L-H$ determine the correction to the leading term $\sqrt{2H/a}$.
Denote Euler's constant by $\gamma_{\rm E}$.

\begin{proposition}\label{prop:race-optimum}
For any sequence with $v\to0$ and $a\to\infty$,
\begin{equation}\label{eq:race-ball-expansion}
 D_n(v)=\sqrt{\frac{\log(1+v)}v}-\frac{1-\gamma_{\rm E}+o(1)}{nv}.
\end{equation}
Over all smoothing scales,
\begin{equation}\label{eq:race-optimum}
 \sup_{v>0}D_n(v)=1-\frac{\sqrt{1-\gamma_{\rm E}}+o(1)}{\sqrt n}.
\end{equation}
\end{proposition}

The two contributions to $1-D_n(v)$ in \eqref{eq:race-ball-expansion} are, to first order, $v/4$ and $(1-\gamma_{\rm E})/(nv)$.
Balancing the two gives the scale $v\asymp n^{-1/2}$.
We first derive the expansion, then show that every nearly optimal scale lies in its range.

\begin{lemma}\label{lem:ball-information}
If $v\to0$ and $a\to\infty$, then
\begin{equation}\label{eq:race-ball-moments}
\begin{aligned}
 \E(L-H) &= \frac12\log(2\pi a)-\frac12+o(1), \\
 \E\bigl((L-H)^2\bigr) &\sim a, \\
 \|L-H\|_p &\le C_p\sqrt a
 \end{aligned}
\end{equation}
for every fixed $1\le p<\infty$ and all sufficiently large $n$.
\end{lemma}
\begin{proof}
Scale lengths by $\sqrt n/r_n$.
Then $Y_B$ becomes $Y=X+\sqrt v Z$, with $X$ uniform in the ball of radius $\sqrt n$ and $Z$ an independent standard normal vector.
Write $R=\lVert Y\rVert$ and $Q(r)=\Prb(\lVert re_1+\sqrt v Z\rVert\le\sqrt n)$, so that $L=-\log Q(R)$.
Use the logarithm of the relative ball volume as radial coordinate,
\[
 U=n\log(R/\sqrt n)
 \quad\text{and}\quad
 W=(U-H)/\sqrt a .
\]
In this coordinate the radial volume element, divided by the volume of the ball of radius $\sqrt n$, is $e^u\,du$.
Thus $U$ has density $e^u Q(\sqrt n e^{u/n})$.
By Pr\'ekopa's theorem and radial symmetry, $\log Q$ is concave and nonincreasing on $[0,\infty)$ \cite{Prekopa1973}.
Its composition with $u\mapsto\sqrt n e^{u/n}$ is therefore concave, so the density of $U$ is log-concave.
The change of variables gives $W$ a positive continuous density $f_n$ and yields the exact identity
\[
 L=H+\sqrt a\,W+\frac12\log a-\log f_n(W).
\]

To find the fluctuations of $U$, let $T=(R^2-n(1+v))/2$.
Write $\lVert X\rVert=\sqrt n\,e^{-E/n}$ with $E$ a unit exponential, and let $G$ be the component of $Z$ in the direction of $X$.
Then
\[
 T=\frac n2(e^{-2E/n}-1)+\sqrt a\,e^{-E/n}G
        +\frac v2(\lVert Z\rVert^2-n).
\]
Here $G$ is standard normal and independent of $E$.
The first term has absolute value at most $E$, and $\Var(\lVert Z\rVert^2)=2n$.
It follows that $\E T=-n/(n+2)$ and $T/\sqrt a-G\to0$ in $L^2$.
On $R\ge\sqrt n/2$, Taylor's formula gives $U-H=T/(1+v)+O(T^2/n)$ with mean absolute remainder $O(v)$.
On $R<\sqrt n/2$, both $|T|$ and $H$ are $O(n)$ and $\Prb(U<-x)\le e^{-x}$ for $x\ge0$, since $Q\le1$.
The contribution from $U<-n\log2$ is therefore exponentially small.
Consequently $\E(U-H)=-1+o(1)$ and $W$ converges in law to a standard normal.

To pass from convergence in law to convergence of the entropy, we use the following one-dimensional argument (see Proposition~2(c) in \cite{CuleSamworth2010} for the general result).
Let $G_n$ be the distribution function of $W$.
Since $f_n$ is log-concave, applying Pr\'ekopa's theorem to $f_n(w)\ind_{\{w\le x\}}$ and integrating in $w$ shows that $G_n$ is log-concave.
At each $x$, the derivative of $\log G_n$ lies between its secant slopes on $[x,x+h]$ and $[x-h,x]$ for every $h>0$.
Since $G_n\to\Phi$ pointwise, letting $n\to\infty$ and then $h\downarrow0$ gives $(\log G_n)'(x)\to(\log\Phi)'(x)$, hence $f_n(x)\to\varphi(x)$.

The secants of $\log f_n$ on $[-2,-1]$ and $[1,2]$ now have slopes tending to $3/2$ and $-3/2$, respectively.
A concave function lies below the extension of each secant outside its defining interval, so these two lines give $f_n(w)\le C e^{-|w|}$ for all $w$ and all large $n$.
Dominated convergence gives $\E W^2\to1$ and bounded fixed moments of $W$.
It also gives
\[
 -\E\log f_n(W)\to\frac12\log(2\pi e)
 \quad\text{and}\quad
 \E|\log f_n(W)|^p=O(1),
\]
since $u|\log u|^p\le C_p\sqrt u$ on the bounded range of $f_n$.
Substitution in the identity for $L$ proves \eqref{eq:race-ball-moments}.
\end{proof}

\begin{proof}[Proof of Proposition~\ref{prop:race-optimum}]
The normal maximum expansion \cite{Nair1981} gives
\[
 m(z)=\sqrt{2z}+
 \frac{\gamma_{\rm E}-1/2\cdot \log(4\pi z)+o(1)}{\sqrt{2z}}
 \quad\text{as}\quad
 z\to\infty.
\]
Stochastic monotonicity extends the formula for integer sample sizes to real $e^z$ by rounding the sample size in both directions.
Let $\Delta=L-H$.
Since $2H\sim a$, the global Taylor estimate $\sqrt{1+x}=1+x/2-x^2/8+O(|x|^3)$ for $x\ge-1$ and Lemma~\ref{lem:ball-information} give
\[
 \E\sqrt{2L}=\sqrt{2H}+
   \frac{\E\Delta-\E(\Delta^2)/(4H)+o(1)}{\sqrt{2H}}.
\]
For the remaining term $m(L)-\sqrt{2L}$ one may replace $L$ by $H$ with mean error $o(a^{-1/2})$.
Indeed, on $|\Delta|\le H/2$ the derivative of the explicit logarithmic correction is $O((\log H)/H^{3/2})$ and $\E|\Delta|=O(\sqrt a)$.
The remainder there is uniformly $o(a^{-1/2})$.
On the complement the correction $m(z)-\sqrt{2z}$ is bounded on $[0,\infty)$ and $\Prb(|\Delta|>H/2)=O(a^{-1})$.
Thus
\[
 \E m(L)=\sqrt{2H}+\frac{1}{\sqrt{2H}}
 \left(\E\Delta-\frac{\E(\Delta^2)}{4H}+\gamma_{\rm E}
              -\frac12\log(4\pi H)+o(1)\right).
\]
The logarithms cancel since $2H/a\to1$.
The mean shift and the quadratic term in \eqref{eq:race-ball-moments} each contribute $-1/2$, leaving $\gamma_{\rm E}-1$.
Division by $\sqrt a$ gives \eqref{eq:race-ball-expansion}.

We next show that every sequence with $D_n(v)\to1$ lies in the range of this expansion.
For $u=\bar\Phi(c)$, Jensen's inequality for the upper normal tail gives $e^{t\J(u)/u}\le e^{t^2/2}/u$ for every $t>0$.
Optimizing in $t$ and using $F\le\J$ gives
\begin{equation}\label{eq:score-upper}
 0\le m(-\log u)\le\J(u)/u\le\sqrt{-2\log u}
\end{equation}
for $0<u\le1$.
A second application of Jensen's inequality yields $D_n(v)^2\le2\E L/(nv)$.
In the original coordinates $\E L$ is the entropy of $Y_B$.
A Gaussian maximizes entropy at fixed covariance, and $\operatorname{Cov}(Y_B)=(r_n^2/n)(n/(n+2)+v) I_n$.
Hence
\[
 D_n(v)^2\le\frac{\log(1+v)+\log(2\pi e r_n^2/n)}v .
\]
Stirling's formula gives $2\pi e r_n^2/n\to1$.
Since $\log(1+v)/v$ decreases strictly from one, $D_n(v)\to1$ forces $v\to0$.

To exclude bounded $a$, consider the radii at which $q$ drops from $1/4$ to $1/8$.
Gaussian isoperimetry applied to $(B-y)/\sigma$, followed by integration in $y$, gives $\int\J(q)\le\sigma\,\area(\partial B)=\sigma b_n$ \cite{Borell1975}, so $D_n(v)\le1$.
For $0<v<1/2$, let $r_0=r_n\sqrt{1-2v}$ and let $Z$ be a standard Gaussian vector.
If $\lVert y\rVert\le r_0$, then $y+\sigma Z\in B$ whenever $\langle y,Z\rangle\le0$ and $\lVert Z\rVert^2\le2n$, since $n\sigma^2=r_n^2v$.
By symmetry and Markov's inequality this event has probability at least $1/4$.
Thus $q\ge1/4$ on the ball of radius $r_0$.
By Lemma~\ref{lem:score-gap}, $c=\min_{1/8\le u\le1/4}(1-F(u)/\J(u))>0$.
On this range of $q$, the gradient estimate from the proof of \eqref{eq:ball-profile} gives $\J(q)-F(q)\ge c\sigma\lVert\nabla q\rVert$.
As $q$ decreases radially from $1/4$ to $1/8$, each sphere has area at least $b_n(r_0/r_n)^{n-1}$, so radial integration gives
\[
 1-D_n(v)\ge\frac{1}{\sigma b_n}\int_{\R^n}\bigl(\J(q(y))-F(q(y))\bigr)\,dy
 \ge\frac c8(1-2v)^{(n-1)/2}.
\]
The right side stays bounded away from zero along any sequence with $a=nv$ bounded, so $D_n(v)\to1$ also forces $a\to\infty$.

Let $d=1-\gamma_{\rm E}>0$.
In the range of \eqref{eq:race-ball-expansion},
\[
 1-D_n(v)=(1+o(1))\frac v4+\frac{d+o(1)}{nv}.
\]
The choice $v=2\sqrt{d/n}$ gives a loss $(\sqrt d+o(1))/\sqrt n$.
Conversely, any sequence within $1/n$ of the supremum in value has $D_n(v)\to1$, and the arithmetic--geometric mean inequality applied to these two positive terms gives the matching lower bound on its loss.
This proves \eqref{eq:race-optimum}.
\end{proof}

Proposition~\ref{prop:race-ball} and \eqref{eq:race-optimum} give, uniformly over all $n$-dimensional parallelohedra,
\[
 \frac{\iq(K)}{\sqrt2\,b_n}\ge
    1-\frac{\sqrt{1-\gamma_{\rm E}}+o(1)}{\sqrt n}.
\]
The first asymptotic correction to the three normalized infima in Theorem~\ref{thm:main} remains undetermined.

\section{Explicit bounds in each dimension}\label{sec:finite-estimates}

\subsection{The lower bound}
For $n\ge2$ and $v>0$ define
\[
 B_n(v)=\sqrt2\left(
  \frac{1-1/n}{\sqrt{n/(n+2)+v}}
  -\frac{2}{\sqrt{nv\bigl(3/2+(n/2)\log(1+v/2)\bigr)}}\right).
\]

\begin{theorem}\label{thm:finite}
Every parallelohedron $K\subset\R^n$ with $n\ge2$ satisfies $\iq(K)/b_n\ge\max\{1,B_n(v)\}$ for every $v>0$.
At $v=\sqrt{8/n}$ this gives
\begin{equation}\label{eq:finite}
 \frac{\iq(K)}{b_n}\ge\max\left\{1,\sqrt2-\frac4{\sqrt n}\right\}.
\end{equation}
\end{theorem}
\begin{proof}
Subtract \eqref{eq:ball-score-loss} from \eqref{eq:ball-profile} and use Proposition~\ref{prop:race-ball} to obtain the bound $B_n(v)$.
The Euclidean isoperimetric inequality gives the lower bound $1$.

For the second assertion take $v=\sqrt{8/n}$ and $a=nv$.
Since $\log(1+v/2)\ge v/2-v^2/8$, we have $(n/2)\log(1+v/2)\ge a/4-1/2$.
The subtracted term in parentheses is therefore at most $4/a=v/2$.
Using $(1+t)^{-1/2}\ge1-t/2$ with $t=v-2/(n+2)$ bounds the first term from below by
\[
 \left(1-\frac1n\right)\left(1-\frac v2+\frac1{n+2}\right)
 =1-\frac v2+\frac1n\left(\frac v2-\frac3{n+2}\right)
 \ge1-\frac v2,
\]
since $n+2\ge2\sqrt{2n}$.
Thus $B_n(v)\ge\sqrt2(1-v)=\sqrt2-4/\sqrt n$.
\end{proof}

Using the full bound $B_n(\sqrt{8/n})$, we obtain $\iq(K)/b_n>1.0589$ in dimension $100$ and $\iq(K)/b_n>1.2929$ in dimension $1000$.

\subsection{The upper bound}
Combining the unit cube with the mean estimate below gives the upper bound
\[
 U_n=\min\left\{\frac{2n}{b_n},\sqrt2+\frac5{\sqrt n}+\frac1{n-1}\right\}.
\]

\begin{proposition}\label{prop:upper}
For every $n\ge2$, the normalized quotient $R_n$ of a Haar-random covolume-one lattice Voronoi cell satisfies
\[
 \E R_n\le\sqrt2+\frac5{\sqrt n}+\frac1{n-1}.
\]
Consequently $\VV_n/b_n\le U_n$.
\end{proposition}
\begin{proof}
The pyramid identity \eqref{eq:cone-identity} gives $\E R_n=\E H^{-1}$.
Integrate $\Prb(H^{-1}>u)$ using the first bound of \eqref{eq:height-tails} from $\sqrt2$ to $2$ and the second from $2$ to infinity.
Since $\int_2^\infty(2/u)^n\,du/2=1/(n-1)$, this gives
\[
 \E R_n\le\sqrt2+
 \int_{\sqrt2}^2\min\left\{1,\frac{17}{4}e^{-n(1-2/u^2)^2/6}\right\}\,du
 +\frac1{n-1}.
\]
For $\sqrt2\le u\le2$ the factor $(u+\sqrt2)/u^2$ is at least $(2+\sqrt2)/4>4/5$, so $1-2/u^2\ge(4/5)(u-\sqrt2)$.
With $z=(4/5)\sqrt n\,(u-\sqrt2)$ the remaining integral is at most
\[
 \frac{5}{4\sqrt n}\int_0^\infty\min\left\{1,\frac{17}{4}e^{-z^2/6}\right\}\,dz
 \le\frac{5}{4\sqrt n}\left(3+\int_3^\infty e^{3-z}\,dz\right)
 =\frac5{\sqrt n}.
\]
The last inequality uses $17/4\le e^{3/2}$ and $z^2/6\ge z-3/2$.
Finally $\VV_n/b_n\le\E R_n$, and the Voronoi cell of $\Z^n$ is the unit cube, with normalized quotient $2n/b_n$.
\end{proof}

\begin{figure}[!h]
\centering
\tikzset{every node/.append style={scale=1.055555556}}
\resizebox{.90\textwidth}{!}{\input{bounds.tex}}
\caption{The lower bound $\max\{1,B_n(\sqrt{8/n})\}$ and upper bound $U_n$ for the normalized isoperimetric quotients at sampled integer dimensions from $10^2$ to $10^{10}$.}
\label{fig:bounds}
\end{figure}

Figure~\ref{fig:bounds} compares the two bounds at $v=\sqrt{8/n}$.
For every $n\ge2$ and $v>0$, the two estimates give
\begin{equation}\label{eq:finite-sandwich}
 \max\{1,B_n(v)\}
 \le\frac{\PP_n}{b_n}\le\frac{\AV_n}{b_n}\le\frac{\VV_n}{b_n}\le U_n.
\end{equation}

\section{A convex integer tile with small projections}\label{sec:integer-projections}

For a convex body $K$ and a unit vector $u$, let $p_K(u)$ be the $(n-1)$-dimensional volume of the orthogonal projection of $K$ onto $u^\perp$.
Write $M(K)=\max_{\lVert u\rVert=1}p_K(u)$ and extend $p_K$ to all vectors by positive homogeneity, with $p_K(0)=0$.

\begin{theorem}\label{thm:integer-projections}
There is an absolute constant $C$ such that in every dimension $n\ge 2$ there is a convex $\Z^n$-tile $K$ satisfying
\[
 M(K)\le C\sqrt{\log n}
 \quad\text{and}\quad
 \area(\partial K)\le C\sqrt{n\log n}.
\]
\end{theorem}

To control all projections, return to the ray through a uniform point of a volume-one Voronoi cell $P=\Vor L$ in Section~\ref{sec:upper}.
Write $h=r_nH$ for the distance from the origin to the plane of the facet hit by the ray.
The pyramid identity \eqref{eq:cone-identity} gives $\area(\partial P)=n\E(h^{-1}\mid L)$.
The second moment gives the directional bound
\begin{equation}\label{eq:projection-height}
 M(P)\le\frac{\sqrt n}{2}\bigl(\E(h^{-2}\mid L)\bigr)^{1/2}.
\end{equation}
Indeed, let $A_F$, $\nu_F$ and $h_F$ denote the area, outward unit normal and height of a facet $F$.
If $v$ is the neighboring lattice point across $F$, reflection $x\mapsto v-x$ preserves $F$, so its centroid is $v/2=h_F\nu_F$.
For each unit vector $u$, the divergence theorem applied to $x\mapsto\langle x,u\rangle u$ now gives
\[
 \sum_F A_Fh_F\langle\nu_F,u\rangle^2=\vol(P)=1.
\]
By Cauchy's projection formula and Cauchy--Schwarz,
\[
 2p_P(u)=\sum_F A_F|\langle\nu_F,u\rangle|
 \le\left(\sum_F A_F/h_F\right)^{1/2}.
\]
The ray hits $F$ with probability $A_Fh_F/n$, which gives \eqref{eq:projection-height}.

We seek a lattice $L=B\Z^n$ for which $M(P)=O(1)$ and the change of variables back to $\Z^n$ costs little.
For an invertible map $T$, the transformation of facet area in Cauchy's formula gives $p_{TK}(u)=|\det T|\,p_K(T^{-1}u)$.
Thus if $\det B=1$, the integer tile $B^{-1}P$ satisfies $M(B^{-1}P)\le\lVert B\rVert_{\mathrm{op}}M(P)$.

For large $n$ let $d=n-1$, $m=\lceil\sqrt n\rceil$, $R=A\lceil\sqrt{\log_2 n}\rceil$ for a sufficiently large fixed integer $A$, and $\delta=R^{-d}$.
Choose the coordinates of $\alpha\in\R^d$ independently and uniformly from $[-1/m,1/m)$.
Let
\[
 B=\begin{pmatrix}R I_d&R\alpha\\ 0&\delta\end{pmatrix}
 \quad\text{and}\quad
 v_{z,j}=B(z,j)=(R(z+j\alpha),\delta j)
\]
for $(z,j)\in\Z^d\times\Z$.
Then $\det B=1$ and $\lVert B\rVert_{\mathrm{op}}\le2R$ since $\lVert\alpha\rVert\le1$.
The vectors $v_{z,j}$ with a given $j$ form one horizontal layer.
For the layers nearest the origin we estimate a Gaussian sum, and for the more distant layers the random horizontal shift gives an integral bound.

\begin{lemma}\label{lem:integer-layers}
Let $J=n^{12}$ and let $T_0$ denote the sum of $e^{-\lVert v_{z,j}\rVert^2/16}$ over $|j|<J$ and $(z,j)\ne(0,0)$.
For $A$ sufficiently large, all sufficiently large $n$ and every nonnegative integrable radially nonincreasing function $f$,
\[
 \E T_0\le n^{-40}
 \quad\text{and}\quad
 \E\sum_{|j|\ge J}\sum_{z\in\Z^d} f(v_{z,j})
 \le (1+C n^{-10})\int_{\R^n} f(x)\,dx .
\]
Both estimates also hold if each $\alpha_i$ is rounded down to a multiple of $2^{-n^2}/m$.
\end{lemma}
\begin{proof}
Let $a=R^2/16$ and $q_a(t)=\sum_{k\in\Z}e^{-a(k+t)^2}$.
For $j\ne0$ let $s=|j|/m$.
Omitting the last coordinate and using independence gives
\[
 \E\sum_{z\in\Z^d} e^{-\lVert v_{z,j}\rVert^2/16}
 \le \left(\frac1{2s}\int_{-s}^s q_a(t)\,dt\right)^d .
\]
Divide the interval of integration into $|j|$ intervals of length $2/m$.
When $m\ge4$, each meets at most one interval of radius $1/(2m)$ about an integer, so at least half its length lies at distance at least $1/(2m)$ from $\Z$.
Keeping a nearest-integer term gives $q_a(t)\le e^{-a\,\operatorname{dist}(t,\Z)^2}+Ce^{-a/4}$ for $a\ge1$, since the remaining distances are bounded below in pairs by $1/2,3/2,\ldots$.
The average is therefore at most $(1+e^{-a/(4m^2)})/2+Ce^{-a/4}\le1-ca/m^2$ for sufficiently large fixed $A$ and all large $n$, since $a\asymp A^2\log n$ and $m^2\asymp n$.
Raising this bound to the power $d=n-1$ bounds the contribution of each nonzero layer by $e^{-c'a}$.
The layer $j=0$ contributes $q_a(0)^d-1=O(ne^{-a})$ after omission of the origin.
Taking $A$ sufficiently large makes the total over fewer than $2J$ layers at most $n^{-40}$ in expectation.

For $|j|\ge J$, the density of $j\alpha_i$ modulo one is at most $1+m/(2|j|)$, by counting integer translates of an interval of length $2|j|/m$.
The joint density on the $d$-dimensional torus is therefore at most $(1+m/(2J))^d\le1+C n^{-10}$.
Let $g(t)=\int_{\R^d}f(x,t)\,dx$.
Integrating over a horizontal period and using $R^{-d}=\delta$ bounds the expected sum in layer $j$ by $(1+C n^{-10})\delta g(j\delta)$.
Since $g$ is even and nonincreasing on $[0,\infty)$, comparison with the integral gives $\delta\sum_{j\ne0}g(j\delta)\le\int_{\R}g(t)\,dt$.

Rounding $\alpha$ to $\alpha'$ sends each lattice vector $x$ to $Tx=x+(R^n(\alpha'-\alpha)x_n,0)$.
Since $\lVert\alpha'-\alpha\rVert\le2^{-n^2}$, we have $\lVert Tx\rVert\ge(1-\varepsilon)\lVert x\rVert$ with $\varepsilon=R^n2^{-n^2}\le n^{-12}$ for large $n$.
The rounded Gaussian sum is bounded by the continuous one with denominator $32$ in place of $16$, for which the same proof applies with $a=R^2/32$.
For the distant layers apply the continuous estimate to $x\mapsto f((1-\varepsilon)x)$, which bounds $f(Tx)$ and changes the integral by a factor $(1-\varepsilon)^{-n}=1+O(n^{-11})$.
\end{proof}

\begin{proof}[Proof of Theorem~\ref{thm:integer-projections}]
Let $V=n^5$ and let $D_V$ be the centered ball of volume $V$.
The radius $b=r_n V^{1/n}$ of $D_V$ satisfies $b\le\sqrt n$ for large $n$, since $r_n/\sqrt n\to1/\sqrt{2\pi e}$.
Choose $X$ uniformly in $P=\Vor L$ and let $N$ count the points of the coset $X+L$ in $D_V$.
Since $X$ is a point of least norm in this coset, $N=0$ exactly when $X\notin D_V$.
As in Section~\ref{sec:upper}, the conditional mean of $N$ is $V$, and $\E(N(N-1)\mid L)$ is the sum over $v\in L\setminus\{0\}$ of $\eta(v)=\vol(D_V\cap(D_V+v))$.
This overlap is radial and nonincreasing with $\int\eta=V^2$.
The two-ball inclusion from Section~\ref{sec:ball} gives $\eta(v)\le V e^{-\lVert v\rVert^2/8}$ since $b^2\le n$.
Applying Lemma~\ref{lem:integer-layers} to the two ranges of $j$ and averaging over both the lattice and the coset gives
\[
 \Var N\le V+V n^{-40}+C n^{-10}V^2\le2V .
\]
Writing $\varepsilon_L=\Prb(X\notin D_V\mid L)$, we obtain $\E\varepsilon_L=\Prb(N=0)\le2/V$.

Let $\lambda=r_n n^{-4/n}$, the radius of the ball of volume $n^{-4}$.
With either choice of coordinates in Lemma~\ref{lem:integer-layers}, we can now choose a lattice for which
\[
 \begin{aligned}
 T_0 &\le n^{-20},\\
 \lVert v_{z,j}\rVert &>\lambda \quad\text{whenever } |j|\ge J,\\
 \varepsilon_L &\le n^{-2}.
 \end{aligned}
\]
The respective failure probabilities are at most $n^{-20}$ by Markov's inequality, $(1+C n^{-10})n^{-4}$ by the estimate for distant layers applied to the ball of radius $\lambda$, and $2n^{-3}$ by Markov's inequality applied to $\varepsilon_L$.
Thus all three properties hold with probability $1-O(n^{-3})$.
The first two properties imply that every nonzero lattice vector has length at least one for large $n$, hence $h\ge1/2$.

A facet of height $h_F<r_n/4$ has label $v$ of length $2h_F<r_n/2<\lambda$ for large $n$ and therefore comes from a layer with $|j|<J$.
Cut the pyramid over this facet by the parallel plane at distance $h_F/\sqrt2=\lVert v\rVert/(2\sqrt2)$ from the origin.
The smaller pyramid at the origin occupies the fraction $2^{-n/2}$ of the volume, so at least half the pyramid lies on the facet side of the plane.
Inside $D_V$ this part lies in a cap whose volume, by Section~\ref{sec:upper}, is at most
\[
 \frac V2\left(1-\frac{\lVert v\rVert^2}{8b^2}\right)^{n/2}
 \le \frac V2 e^{-\lVert v\rVert^2/16}.
\]
The pyramids have disjoint interiors, and their parts outside $D_V$ have total volume at most $\varepsilon_L$.
Summing the cap volumes and doubling gives $\Prb(h<r_n/4\mid L)\le 2\varepsilon_L+VT_0$.
Since $h\ge1/2$,
\[
 \E(h^{-2}\mid L)\le \frac{16}{r_n^2}+4(2n^{-2}+Vn^{-20})=O(1/n).
\]
Equation~\eqref{eq:projection-height} gives $M(P)=O(1)$, so $K=B^{-1}P$ tiles by $\Z^n$ and satisfies $M(K)\le 2R M(P)=O(\sqrt{\log n})$.
Averaging Cauchy's projection formula over unit directions gives mean projection volume $\omega_{n-1}\area(\partial K)/(n\omega_n)$, with $\omega_{n-1}/(n\omega_n)\sim1/\sqrt{2\pi n}$.
This proves the asserted surface area bound for the same tile.
The unit cube covers the finitely many remaining dimensions $n\ge 2$ after enlarging $C$.
\end{proof}

\begin{corollary}\label{cor:integer-oracle}
In every dimension $n\ge 2$, a randomized polynomial-time algorithm samples a nonsingular rational matrix $B$ of bit length polynomial in $n$.
Let $K=B^{-1}\Vor(B\Z^n)$.
With probability $1-O(n^{-3})$ over the choice of $B$, the tile $K$ satisfies both bounds in Theorem~\ref{thm:integer-projections} and the following optimization guarantee holds simultaneously for every $u\in\mathbb Q^n$.

Given such a matrix $B$ and any $u\in\mathbb Q^n$, a deterministic algorithm returns a point $x\in K\cap\mathbb Q^n$ in time polynomial in $n$ and the bit length of $u$, with
\[
 \langle u,x\rangle\ge
 \left(\frac14-C\frac{\log n}{n}\right)
 \max_{y\in K}\langle u,y\rangle
\]
for an absolute constant $C$.
For each fixed $0<c<1/4$, the approximation factor can be taken to be $c$ in every dimension $n\ge 2$.
\end{corollary}

\begin{proof}
Use the rounded coordinates in Lemma~\ref{lem:integer-layers}.
Each is uniform on $2^{n^2+1}$ rational grid points, so $n^2+1$ fair bits per coordinate suffice to sample $B$ in polynomial time.
Its bit length is polynomial since $R$ is an integer and $\delta=R^{1-n}$.
Let $n$ be large and suppose that the three properties in the proof of Theorem~\ref{thm:integer-projections} hold for $L=B\Z^n$.

We turn the estimate $\varepsilon_L\le n^{-2}$ into a bound on the covering radius using the doubling argument in Lemma~2.5 of~\cite{OrdentlichRegevWeiss2022}.
Since $X$ represents a uniform coset, the image of $D_V$ in $\R^n/L$ has measure at least $1-n^{-2}>1/2$ and intersects every translate of itself.
Taking differences, every coset has a representative in $D_V-D_V=2D_V$.
The points of $P=\Vor L$ are nearest representatives, so $P\subseteq2b B_1^n$ with $b=r_n V^{1/n}$.

Let $a=\lambda/2$ with $\lambda=r_n n^{-4/n}$ as in the preceding proof.
Every vector in a layer $|j|\ge J$ has length greater than $2a$, so its Voronoi inequality holds throughout $a B_1^n$.
Truncation at this radius would retain the fraction $a/(2b)$ of $P$ in every direction.
We use a rational polytope just inside this ball to optimize by linear programming.

Choose a rational $s$ with $a^2-2\le s\le a^2-1$.
To compute it in polynomial time, use $\omega_0=1$, $\omega_1=2$ and $\omega_k=2\pi\omega_{k-2}/k$ with rational bounds on $\pi$ differing by at most $n^{-2}$, obtainable from the alternating series for $4\arctan 1$.
The resulting bounds for $\omega_n$ have ratio $1+O(n^{-1})$, so substituting them in $\omega_n^2(4t)^n=n^{-8}$ brackets the root $a^2\asymp n$ to within $O(n^{-1})$.
Rational bisection to fixed absolute accuracy then suffices.
Define
\[
 g(t)=\max_{k\in\Z}(2kt-k^2)
 \quad\text{and}\quad
 D=\Bigl\{y\in\R^n:\sum_{i=1}^n g(ny_i)\le n^2s\Bigr\}.
\]
Since $t^2-g(t)$ is the squared distance from $t$ to the nearest integer, it lies between $0$ and $1/4$.
Thus $\sqrt s\,B_1^n\subseteq D\subseteq\sqrt{s+1/(4n)}\,B_1^n\subseteq aB_1^n$.
In particular
\[
 \frac{\sqrt s}{2b}P\subseteq P\cap D\subseteq P
 \quad\text{with}\quad
 \frac{\sqrt s}{2b}=\frac{1-O(n^{-1})}{4n^{9/n}}
 \ge\frac14-C\frac{\log n}n.
\]

Only the $2J-1$ layers with $|j|<J$ can cut $D$.
In layer $j$, write $w_{ijk}=R(k+j\alpha_i)$ for the possible values of coordinate $i\le d$.
These coordinates can be chosen independently, so all the Voronoi inequalities from the layer are equivalent to
\[
 \sum_{i=1}^d\max_{k\in\Z}\bigl(2w_{ijk}y_i-w_{ijk}^2\bigr)
       +2\delta j y_n-(\delta j)^2\le0.
\]
The maximizing $k$ is a nearest integer to $y_i/R-j\alpha_i$.
Since $a<n$ for large $n$, we may impose $|y_i|\le n$ and take only $|k+j\alpha_i|\le n+1$ in these maxima and $|k|\le n^2$ in $g(ny_i)$.
Replace each maximum by an auxiliary variable bounded below by its affine terms.
The remaining inequalities bound sums of these variables from above, so the resulting rational linear system projects exactly onto $P\cap D$.
It has polynomial size and can be written down in polynomial time from $B$.
Maximizing $\langle B^{-T}u,y\rangle$ gives a rational optimizer in polynomial time \cite{Khachiyan1980}, and $x=B^{-1}y$ has the claimed support value.
For each fixed $0<c<1/4$, take $B=I_n$ in the finitely many remaining dimensions $n\ge 2$ and optimize exactly over the centered unit cube.
Enlarging the constants in Theorem~\ref{thm:integer-projections} gives the geometric bounds in these dimensions.
\end{proof}

We end with two questions.
Naor and Regev call a convex body whose integer translates tile $\R^n$ an \emph{integer parallelotope}.
The unit cube $[0,1]^n$ is one example and has surface area $2n$.
In Question 2 of \cite{NaorRegev2023}, they ask whether one can find such a tile in every dimension with surface area of order $\sqrt n$, the same order as that of a ball of volume one.
We ask whether the integer lattice also permits the sharp factor from Theorem~\ref{thm:main}.

\begin{question}\label{q:integer-sharp-constant}
Does there exist a convex $\Z^n$-tile in each dimension $n$ such that the ratio of its surface area to that of a ball of equal volume tends to $\sqrt2$ as $n\to\infty$?
\end{question}

We measure the isoperimetric deficit of a tile by the ratio of its surface area to that of a ball of equal volume.
Once the tiling rule is fixed, an \emph{optimal tile} is one with the smallest deficit in its dimension.
For convex lattice tiles, the known minima in dimensions two, three and four give successively larger values.
Figure~\ref{fig:bounds} shows how the lower and upper bounds close in on $\sqrt2$ by plotting their numerical values at selected dimensions. How the minimum deficit changes from one dimension to the next remains open.
Milman suggested to us the following natural question.

\begin{question}[Milman]\label{q:dimension-monotonicity}
Is the isoperimetric deficit of an optimal tile monotone with dimension?
\end{question}

\section*{Acknowledgments}

\subsection*{Human interactions}
The authors wish to thank Annalisa Cesaroni for stimulating discussions on a preliminary version of this paper. 
We are grateful to Emanuel Milman for helpful discussions which sparked a joint project on a generalization of this work.
We also thank Frank Morgan for useful suggestions.

\subsection*{Machine interactions}
The authors wish to give LLMs the credit they deserve.

The sharp asymptotic constant $\sqrt2$ was proved by the authors entirely by hand. 
After ChatGPT 6.2 Astra was released, we asked it if it could improve our proof.

It significantly improved the lower bound error to $O(n^{-1/2})$ by combining an exponential clock comparison and rearrangement to a ball with the fact that two balls of radius $r$ whose centers are distance $d\le 2r$ apart intersect inside a ball centered at their midpoint with radius $\sqrt{r^2-d^2/4}$.

\vspace{-4pt}
It improved the upper bound error from $O(\sqrt{\log n/n})$ to $O(n^{-1/2})$ by using spherical cap volume estimates to control the probability of hitting low facets \eqref{eq:upper-rate}, \eqref{eq:height-tails}.

It also recasts the upper bound proof in terms of the first moment of inverse facet height \eqref{eq:cone-identity}. The corresponding second moment, together with the centroid property of Voronoi facets, bounds all hyperplane projections \eqref{eq:projection-height}, making the route visible to resolving the question of Naor and Regev \cite{NaorRegev2023}.

ChatGPT 6.2 Astra also helped draw Figures \ref{fig:optimal-shapes} and \ref{fig:bounds}. Claude Fable 5.2 gave us very helpful referee comments.

We would not have been able to present our work with simple and clean proofs around a unified conceptual idea this such without the help of LLMs.

\subsection*{Funding and support}
MN is member of INdAM/GNAMPA 
(National Group for Mathematical Analysis, Probability and their Applications, Italy).

\printbibliography

\end{document}